\documentclass[leqno]{prims}
\usepackage{hyperref}
\usepackage{amssymb,amsmath}
\usepackage{enumerate}
\usepackage{graphicx}

\newtheorem{thm}{Theorem}[section]

\newtheorem{lem}[thm]{Lemma}

\theoremstyle{definition}

\numberwithin{equation}{section}  

\newcommand{\N}{\mathbb{N}}
\newcommand{\R}{\mathbb{R}}
\newcommand{\dif}{\mathrm{d}}

\newcommand{\dz}{\dif z}

\newcommand{\diverg}{\mathop{\mathrm{div}}}

\newcommand{\tanpart}[1]{\left[ #1 \right]_{\textrm{tan}}}

\newcommand{\ue}{u^\varepsilon}
\newcommand{\uinit}{{u_{*}}}
\newcommand{\Rem}{{r^\varepsilon}}
\newcommand{\BorderNav}{\partial\Omega \setminus \Sigma}

\newcommand{\Hspace}{L^2_\sigma(\Omega)}

\newcommand{\Omext}{\mathcal{O}}

\begin{document}

\TitleHead{On the controllability of the Navier-Stokes equation}
\title{On the controllability of the Navier-Stokes equation in spite of  boundary layers}

\AuthorHead{Coron, Marbach and Sueur}
\author{Jean-Michel \textsc{Coron}, 
Fr\'ed\'eric \textsc{Marbach}
and
Franck \textsc{Sueur}}



\maketitle

\begin{abstract}

In this proceeding we expose a particular case of a  recent result obtained in   \cite{CMS} by the authors regarding the incompressible Navier-Stokes equations in a smooth bounded and simply connected bounded domain, either in 2D or in 3D, with a Navier slip-with-friction boundary condition except on a part of the boundary. 
This under-determination encodes that one has control over the remaining part of the boundary. 
We prove that for any initial data, for any positive time, there exists a weak Leray solution which vanishes at this given time.

\end{abstract}


\section{Geometric setting}

We consider a smooth bounded and simply connected\footnote{Indeed our analysis  also covers the case of a multiply connected domain for some  controls located on a  part of the boundary intersecting all its connected  components, but we will stick here to this simple case and we refer to  \cite{CMS}  for the general case.} domain~$\Omega$ in $\R^d$, with  $d = 2$ or $d = 3$.
 Inside this  domain, an incompressible viscous fluid evolves under the Navier-Stokes 
equations. We will name~$u$ its velocity field and $p$ the associated pressure. 
The equations read:
\begin{equation}
\label{NS-eq}
         \partial_t u + (u \cdot \nabla) u - \Delta u + \nabla p  = 0  
         \quad \textrm{and} \quad  \diverg u  = 0 
         \quad \textrm{in } \Omega .
\end{equation}
Let us emphasize that the fluid density and the viscosity coefficient are set equal to one for the sake of clarity. 

\section{Boundary conditions}
  \label{sec-bd}

For an impermeable wall, it is natural to prescribe the condition $ u \cdot n  = 0$ on~$\partial \Omega$, where  $n$ denotes the outward pointing normal to the domain, 
which means that the fluid cannot escape the domain and that there is no cavitation at the boundary. 
Indeed in the case of a perfect fluid, driven by the Euler equations rather than the Navier-Stokes, such a condition is sufficient to have existence and uniqueness to the Cauchy problem in various appropriate functional settings. 
For the case here of the Navier-Stokes equations, an extra condition has to be added. 
The two following propositions are the most used (in complement to the previous condition) : 
\begin{itemize}
\item the no-slip condition $\tanpart{u} = 0$ (dating
back to Stokes in 1851), where
\begin{equation}
  \label{not-tan}
\tanpart{u}:= u - (u \cdot n) n 
\end{equation}
 denotes the tangential part of the vector field $u$.
\item the slip-with-friction  condition  $N(u) = 0$, where
\begin{equation}
  \label{defN}
N(u) := \tanpart{D(u)n +  \alpha u}
\quad \text{with} \quad
D(u) := \left( \frac{1}{2}\left(\partial_i u_j + \partial_j u_i\right) \right)_{1 \leqslant i,j  \leqslant d} \end{equation}
 the rate of strain tensor (or shear stress) and $\alpha $ is a real constant coefficient for simplicity\footnote{Our analysis  also covers the case where 
 $\alpha$ is a smooth matrix-valued function.}.
This condition dates back to Navier in 1833 (see \cite{navier1823memoire}).
 This coefficient describes the friction near the boundary. 
 Let us observe that, 
 formally, when 
 $\alpha \rightarrow +\infty$, the Navier condition 
 reduces to the usual no-slip  condition.
\end{itemize}
%

\section{The Cauchy problem}

Let us recall the following result, where  $\Hspace$ denotes the closure in $L^2(\Omega)$ of smooth 
divergence free vector fields which are tangent to~$\partial \Omega$. 
\begin{thm}
\label{th-leray}
 Let  $u_0 \in \Hspace$.
 Then there exists a global weak solution $u$ associated with the initial data $u_0$. 
\end{thm}
This result dates back to the pioneering work \cite{Leray} by Leray where it is proved that 
   $u \in \mathcal{C}_w^0([0+\infty);\Hspace) \cap L^2((0,+\infty); H^1(\Omega))$.
   Moreover, Leray proved the following partial regularity property: for almost every $t$ in $(0,+\infty)$, 
   $u(t,\cdot)$ is $\mathcal{C}^\infty (\Omega)$.
    
Even though Leray's paper tackled the case of the no-slip condition, this result  can be adapted almost right away to the case of the Navier slip-with-friction condition (see~\cite[Section 3]{MR2754340}).
   
\section{The control problem}

\begin{figure}[ht!]
   \begin{center}
       \includegraphics{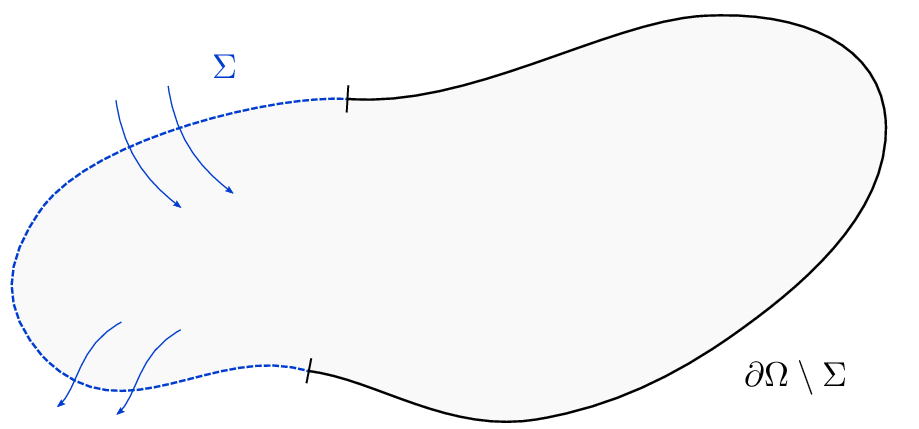}
   \end{center}
   \caption{\label{Figure:control_setting} The control problem}
\end{figure}

We now assume that we are able to act on a non-empty open part  $\Sigma$ of the full boundary~$\partial\Omega$. 
In particular we may let some fluid enter into the domain and the same volume of fluid go out of the domain (recall that the fluid is incompressible).  
Then the setting we have in mind now is the following (see Figure~\ref{Figure:control_setting}).
\begin{itemize}
\item On the part  $\BorderNav$, some boundary conditions are prescribed, either the no-slip condition  $ u = 0$  or the Navier condition $ u \cdot n  = 0$ and $ N(u) = 0$ (that is without any source terms or ability to modify the slip coefficient $\alpha$ which is assumed to be given once and for all).
\item On the part  $\Sigma$, we are free to choose a boundary condition which is relevant for some purpose.
\end{itemize}

More precisely we have in mind to drive the system from an arbitrary initial data to some given state at some given time. 
The following goal, first suggested by Jacques-Louis Lions in the late 80's (cf. for instance \cite{MR1147191}) tackles the case where the target is the rest state. 

\smallskip

\label{lbOP}{\textbf{Open Problem (OP)}}. 
For any $T > 0$ and $u_0$ in $L^2_\sigma (\Omega)$, does there exist a solution to the Navier-Stokes  system with $u(0, \cdot) = u_0$ 
 such that $u(T,\cdot) = 0$ ?
\smallskip

Above the Navier-Stokes system to which \hyperref[lbOP]{(OP)} refers is  constituted of  the incompressible Navier-Stokes  equations \eqref{NS-eq}  in $ \Omega$ and 
of the no-slip condition or the Navier condition on $\BorderNav$,  but without any boundary condition prescribed on the controlled part $\Sigma$ of the boundary. 
Such a system is therefore under-determined so that uniqueness of a solution is not expected (even in the 2D  case for which uniqueness of Leray solutions is known in the uncontrolled setting corresponding to the case where $\Sigma = \emptyset$). 
Indeed in the  formulation above the control is implicit:  a relevant condition to prescribe as a control on  $\Sigma $ can be recovered  by taking the trace on $\Sigma$ of a convenient solution to the under-determined system. 

Observe that there is no restriction regarding the sizes neither of the time $T > 0$ nor of the initial data $u_0$ in $L^2_\sigma (\Omega)$.
In the terminology of control theory a positive result to this question amounts to the small-time global controllability of the 
Navier-Stokes, or more precisely  the   \textit{small-time global exact null controllability}  since the target in \hyperref[lbOP]{(OP)} is the  rest state and has to be reached exactly.

\section{Our result}

 In Lions' original question, the 
boundary condition on the uncontrolled part $\BorderNav$ of the boundary is the no-slip  
boundary condition. Our goal here is to present the following result establishing a positive answer to \hyperref[lbOP]{(OP)} in the case where some Navier conditions 
are prescribed on $\BorderNav$.

\begin{thm} 
\label{main}
Let $T > 0$ and $u_0 \in \Hspace$. There exists a weak solution $u$ 
to 
 \begin{equation} \label{eq.main}
    \left\{
    \begin{aligned}
        \partial_t u + (u \cdot \nabla) u - \Delta u + \nabla p & = 0
        && \quad \textrm{in } \Omega, \\
        \diverg u & = 0
        && \quad \textrm{in } \Omega, \\
        u \cdot n & = 0
        && \quad \textrm{on } \BorderNav, \\
        N(u) & = 0
        && \quad \textrm{on } \BorderNav
    \end{aligned}
    \right.
\end{equation}
satisfying  $u(0,\cdot) = u_0$ and $u(T, \cdot) =  0$.
\end{thm}

Theorem \ref{main} does not require any condition on the coefficient $\alpha$ appearing in the definition \eqref{defN} of $N$. 
Indeed, observe that there is no asymptotic parameter in the statement above.
Still the next lines about the proof of it will be full of  $\varepsilon $.

Let us also mention that the results in  \cite{CMS}  are more general, in particular they prove that one may intercept at time $T$ any  smooth uncontrolled solution to the Navier-Stokes system with Navier condition on the full boundary $\partial \Omega$. 

\section{Earlier results}

When  Jacques-Louis Lions formulated it in the late 80's, \hyperref[lbOP]{(OP)} was pretty impressive
 since the answer  was not even known in the case of  the heat equation. 
For this equation the first key breakthroughs were obtained by 
\cite{LebeauRobbiano}
and \cite{FursikovImanuvilov} thanks to 
 Carleman inequalities respectively  associated with  parabolic and  elliptic second order operators.
The latter has been then extended to the Stokes equations 
and later on to the Navier-Stokes equations in the case of small initial data 
by  Imanuvilov  in~\cite{MR1804497}. The smallness assumption implies that the quadratic convective term may be seen as a perturbation term so that 
the result can be obtained from the controllability of the Stokes equations by a fixed point strategy. 
This result has since been improved in~\cite{MR2103189} by Fern\'andez-Cara, 
Guerrero, Imanuvilov and Puel. 

All these works deal with the case of the no-slip  boundary condition. 
For Navier slip-with-friction boundary conditions, let us mention 
\cite{MR2224824} and \cite{MR2268275} which prove in particular local
null controllability when the initial data is small.

The case of large initial data was first tackled in~\cite{MR1393067}, where the 
first author proves a small-time global result in a 2D setting with Navier
boundary conditions: the smallness obtained within the inside
of the domain is good, but the estimates up to the boundary are not sufficient
 to conclude using a known local result. In fact, when there is no boundary,
the first author and Fursikov prove in~\cite{MR1470445} a {small-time} 
global exact null controllability result
(in this setting, the control is a source term located in a small subset of the domain) 
thanks  to  the return method and to  the global controllability of the incompressible 
Euler equations  (for large smooth initial data). 
Likewise, in~\cite{MR1728643}, Fursikov and Imanuvilov prove 
{small-time} global exact null controllability when the control is supported 
on the whole boundary (i.e. $\Gamma = \partial \Omega$). 
In~\cite{MR2560050}, Chapouly obtains global exact null controllability for 
Navier-Stokes in a 2D rectangular domain
 under Lions' boundary condition (corresponding to the case where $\alpha = 0$ in the Navier condition)
on the uncontrolled part of boundary.  

Still the approaches used in the aforementioned papers 
 failed to deal with the viscous boundary layers  appearing near the uncontrolled part of the boundary.
This is precisely the goal of this paper to promote the 
  \textit{well-prepared dissipation method} in order to  obtain some controllability results despite the presence of boundary layers. 
This method was first introduced in~\cite{MR3227326} by the second author in order to deduce a controllability result for the 1D  Burgers equation.  
The extension of this method to the Navier-Stokes equations will be crucial in our proof of Theorem~\ref{main}. 
 In particular here the method will be implemented thanks to a multi-scale expansion describing the boundary layer occurring in the vanishing viscosity limit of the Navier-Stokes equations. The  application of the method is presented in Section \ref{wpdm}.

\section{A few words of caution} 

Next sections are devoted to the scheme of proof of Theorem \ref{main}. 
We will try to highlight a few key ingredients whereas some technical difficulties will be omitted on purpose for the sake of clarity.
 We refer to \cite{CMS} for a complete proof. 
  
  Let us also mention here that we are not going to really use a control all the time in the sense that it will be relevant on some time intervals to choose as boundary condition on 
  $\Sigma$ the same Navier condition than on $\BorderNav$  so that the system then coincides with the uncontrolled one for which  $\Sigma = \emptyset$.

\section{Reduction to an approximate controllability problem from a smooth initial data} 
\label{sec-red}

\begin{figure}[ht!]
    \begin{center}
        \includegraphics{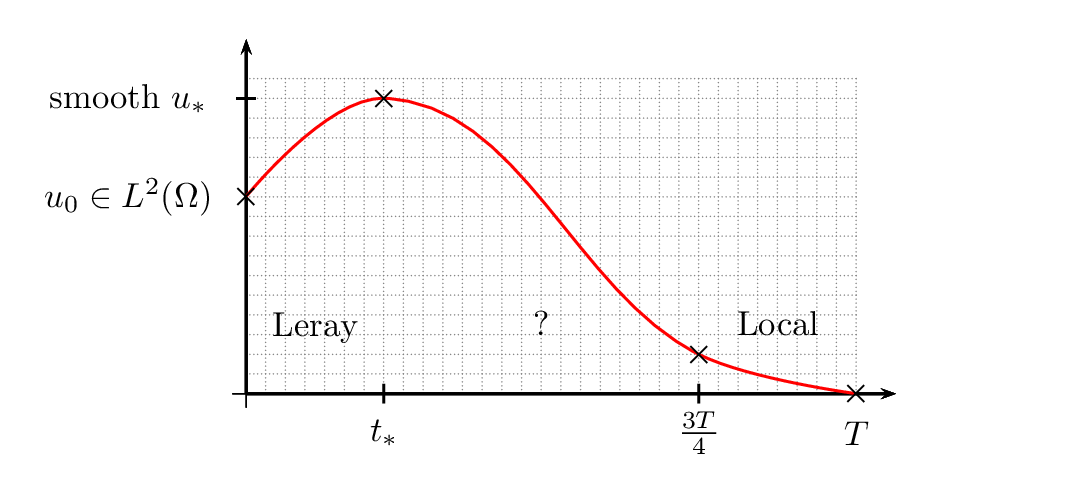}
    \end{center}
    \caption{\label{Figure:graph1} Reduction to a global approximate controllability problem}
\end{figure}

In this section we are going to prove that it is sufficient to have the 
existence of a solution starting from an arbitrary smooth initial data and 
reaching a state close to zero in $L^2_\sigma(\Omega)$, in any positive time in order to conclude 
the proof of Theorem~\ref{main}. 
Indeed according to Leray's partial regularity result hinted above (cf. below 
Theorem~\ref{th-leray}), there exists $t^*$ in $(0,T/2)$ such that $u(t^*,\cdot)$ 
is $\mathcal{C}^\infty (\Omega)$. Let us assume that we are able to prove the 
existence of a solution starting from $u(t^*,\cdot)$ at time $t^*$ and reaching, 
say at time $3T/4$, a state close enough to zero in $L^2_\sigma(\Omega)$ such that the local 
controllability results mentioned above can be applied\footnote{Results available 
in the local controllability literature require to start with an initial data 
which is more regular than $L^2$ but Leray's partial regularization of the 
uncontrolled  Navier-Stokes equations can be used again in order to glue 
these steps. Here we have to pay attention to the preservation of the 
smallness assumption in this regularization argument (cf.~\cite{CMS} for more).}
on the remaining time interval $(3T/4 , T)$. Then the 
concatenation of this three steps yields Theorem~\ref{main} (see Figure~\ref{Figure:graph1}). 
Our task is therefore only to obtain approximate null controllability 
from a smooth initial data on the intermediate time interval $(t^* , 3T/4)$. 
In order to simplify the notations let us pretend that this interval is $(0,T)$ in the next sections. 
On the other hand we will denote $ \uinit$ the initial data, which is 
smooth, for this new problem, in order to distinguish it from the original 
initial data $u_0$ which was only assumed to be in~$L^2_\sigma (\Omega)$. 

\section{A fast and furious control}
\label{sec-ff}

In order to take profit of the nonlinearity at our advantage we aim at reaching approximatively zero thanks to 
 a control which is fast and furious in the sense that its amplitude and duration are scaled with respect to 
a small positive parameter  $\varepsilon$ which is introduced by force and will be ultimately taken small enough. 
Indeed we look for a solution $u$ to \eqref{eq.main} of the form 
\begin{equation}
  \label{scaling}
 u(t, x) =  \frac{ 1 }{\varepsilon }  \ue \left( \frac{ t }{\varepsilon } 
 ,x\right) ,
\end{equation}
having in mind to look for a family of  functions $\ue$ with, typically, 
variations of order $O(1)$ on time interval of order $O(1)$. 
This means for the original  searched solution $u$ having fast transitions on  time interval of order $O(\frac{ 1 }{\varepsilon })$ 
 with furious amplitudes  of order $O(\frac{ 1 }{\varepsilon } )$.
The underlying idea is to start with the ambitious idea to try
to control the system during the shorter time interval $[0, \varepsilon T]$ 
with forcing the system to evolve in a high Reynolds regime. 

 Regarding the pressure $p$ associated with the original solution $u$,
 having in mind Bernoulli's principle which associates the pressure with the square of the velocity, 
  we look for an ansatz of the form 
 $$p(t, x) =   \frac{ 1 }{\varepsilon^2 } p^\varepsilon  \left( \frac{ t }{\varepsilon^2}  ,x \right)   .$$
 This translates the original system in a new system  with 
 four main changes: 
 \begin{enumerate}[i)]
\item The new unknowns  $\ue$ and $p^\varepsilon$ satisfy the Navier-Stokes equations with a small viscosity coefficient $\varepsilon$.
\item The initial data for $\ue$  is now small equal to $\varepsilon \uinit$.
\item The time interval is now $(0, \frac{T}{\varepsilon})$ so that we have to investigate the large time behaviour of the system. In particular, although the initial data is small, nonlinearities will matter. 

\smallskip
The system for $(\ue,p^\varepsilon)$ therefore reads:
\begin{gather}
\label{niou}
   \left\{
    \begin{aligned}
		\partial_t \ue + \left( \ue \cdot \nabla \right) \ue - \varepsilon \Delta \ue + \nabla p^\varepsilon  & =  0  && \quad  \textrm{in}   \quad   (0, T/\varepsilon) \times \Omega,  \\
		 \diverg \ue & =  0 &&  \quad  \textrm{in}  \quad  (0, {T/}{\varepsilon}) \times \Omega,  \\
		\ue\cdot n  & =  0 &&  \quad  \textrm{on} \quad (0, {T/}{\varepsilon}) \times \BorderNav, \\
		 N( \ue )  & =  0  &&  \quad  \textrm{on} \quad  (0, {T/\varepsilon}) \times \BorderNav, \\
		\ue \rvert_{t = 0}  & =  \varepsilon \uinit &&  \quad   \textrm{in} \quad    \Omega.
		 \end{aligned}
    \right.
\end{gather}

\item  Last change but not least, the rest state  is now targeted (at the final time $T/{\varepsilon}$) with more precision. 
Indeed we now plan to prove that there exists a solution $\ue$  to the underdetermined system   \eqref{niou} such that 
\begin{equation}
  \label{newgoal}
 \left\|\ue\left( \frac{T}{\varepsilon} ,\cdot\right)\right\|_{L^2(\Omega)}=o(\varepsilon)
\end{equation}
 in order to deduce from   \eqref{scaling} that there exists a solution $u$ to \eqref{eq.main} such that 
 $\| u( {T},\cdot)\|_{L^2(\Omega)}=o(1)$. In particular, choosing $\varepsilon$  small enough allows to  reach a state arbitrarily close to $0$ in $L^2$.
 This will provide the approximate controllability result mentioned in the previous section.  
\end{enumerate}
%

\section{Inviscid flushing}

When $\varepsilon$ is small,  it is  expected that the analysis of the system   \eqref{niou}  may be built on 
the small-time  global exact controllability of Euler equations.
We therefore consider the following  counterpart of the system   \eqref{niou} 
 where the viscosity term has been dropped out: 
\begin{gather}
  \label{Eeq}
    \left\{
    \begin{aligned}
     \partial_t u^E + \left( u^E \cdot \nabla \right) u^E + \nabla p^E
			 & =   0  &&\quad  \textrm{ in } \Omega,
	\\		  \diverg u^E  & =  0 &&\quad  \textrm{ in } \Omega,
\\   u^E \cdot n  & =  0  &&\quad \textrm{on }  \BorderNav ,
\\                  	u^E \rvert_{t = 0}  & =  \varepsilon \uinit &&\quad   \textrm{in} \quad    \Omega.	
		 \end{aligned}
    \right.
\end{gather}
As the initial data is of order  $O(\varepsilon)$ in $L^\infty$ it is natural to for  a solution  $ u^E$ to  \eqref{Eeq} which is, at least for  times of order $O(1)$,  of the form:
\begin{equation} 
\label{naive}
 u^E  =  \varepsilon u^1 + o( \varepsilon)  \quad   \text{ and }  \quad  p^E = \varepsilon {p^1} + o( \varepsilon)  .
\end{equation}
Plugging expansions \eqref{naive} into \eqref{Eeq}  and grouping terms of order $O_{L^\infty}(\varepsilon)$ 
yields:
\begin{equation}
\label{bougepas}
   \left\{
    \begin{aligned}
\partial_t u^1 + \nabla p^1 &= 0 &&\quad  \textrm{in } \Omega,
\\ \diverg u^1  &= 0 &&\quad \textrm{in } \Omega, 
\\ u^1 \cdot n  &= 0 &&  \quad \textrm{on }  \BorderNav , 
\\  u^1 \rvert_{t = 0}  &= \uinit &&  \quad  \textrm{in } \Omega.
	 \end{aligned}
    \right.
\end{equation}
By elementary combinations of the equations we observe that the system  \eqref{bougepas} does not admit any solution reaching exactly $0$ unless the initial data $ \uinit $ is the gradient of a harmonic function, which is not the case in general. 
System  \eqref{bougepas}  suffers from a lack of controllability which will prevent from using it for our purposes. 

\begin{figure}[ht!]
   \begin{center}
       \includegraphics{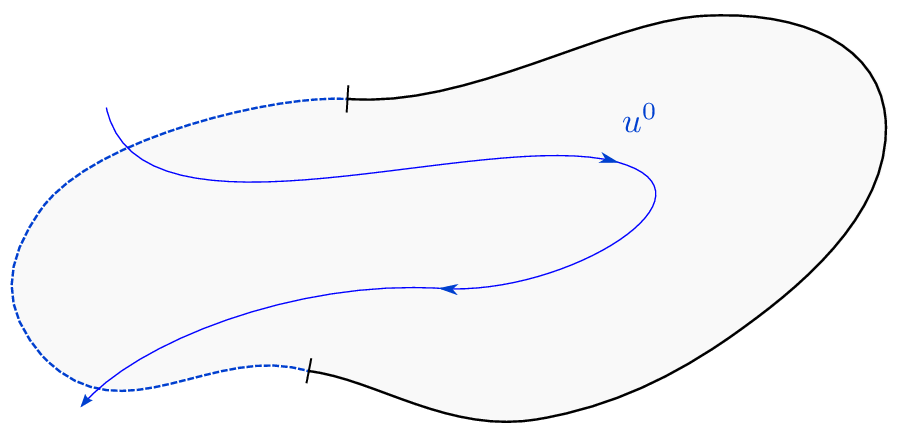}
   \end{center}
   \caption{\label{Figure:control_euler} Inviscid flushing}
\end{figure}

In order to overcome this difficulty we are going to use  the return method first introduced by the first author in \cite{MR1233425}.
 This method takes profit of the nonlinearity thanks to an auxiliary controlled solution to the Euler system.
 Indeed, instead expansions \eqref{naive}, we will rather look for some asymptotic expansions of the form: 
\begin{equation} 
\label{-naive}
 u^E  = u^0 + \varepsilon u^1 +  o( \varepsilon)     \quad\text{and}\quad  
 p^E  = p^0 + \varepsilon p^1 +  o( \varepsilon) , 
\end{equation}
where the extra-term $(u^0,p^0)$ is introduced in order to help to control $(u^1 , p^1 )$. 
Of course $(u^0,p^0)$ has to be solution to the Euler system in order to cancel out  the terms of order $O(1)$  which appear when plugging the expansions 
 \eqref{-naive} into the first three equations of \eqref{Eeq}. 
 Moreover the last equation yields the initial data  $u^0 \rvert_{t = 0}  = 0$  in $ \Omega$. 
 The interest is that the equations obtained by  gathering the terms of order $O(\varepsilon)$  are now: 
\begin{equation}
\label{linu}
  \left\{
    \begin{aligned}
\partial_t u^1 + \left( u^0 \cdot \nabla \right) u^1 +  \left( u^1 \cdot \nabla \right) u^0 + \nabla p^1  &= 0 &&  \quad  \textrm{in } \Omega,
\\  \diverg u^1   &=  0 &&  \quad  \textrm{in } \Omega, 
\\  u^1 \cdot n  & = 0  &&  \quad \textrm{on }  \BorderNav ,
\\  u^1 \rvert_{t = 0}  &= \uinit &&  \quad  \textrm{in } \Omega.
   \end{aligned}
    \right.
\end{equation}
This is the linearisation of the Euler equations around $u^0$ (rather than around $0$ like in \eqref{bougepas}). 
We may now rely on the transport  by $u^0$ (see Figure~\ref{Figure:control_euler}) in order to drive $u^1$ from $\uinit$ to $0$, see the second term in the first equation of \eqref{linu}. 
More precisely we want to use the transport by $u^0$  in order to flush $u^1$ out of the domain.  
Of course the system \eqref{linu} has, in addition to the transport aspect, non local features due to the incompressibility condition.
Still, reasoning on the vorticity of $u^1$, we obtain that if the fluid particles are flushed outside of the physical domain within a time interval of order $O(1)$, say 
$[0,T]$ (that is if any fluid particle initially at $x$ in $\Omega$ moves with the flow  associated with~$u^0$ up to some time $t_x \in (0,T)$ for which it reaches $\Sigma$ with a positive velocity), then $u^1$ can be set equal to $0$.  Observe that this requires a time interval far smaller than the allotted one which is  $[0,T/\varepsilon]$.

On the other hand this auxiliary field $u^0$  also has to vanish at the final time. 
In order to construct such a field, a crucial observation is that the potential flows as solutions to the Euler system enjoy a lot of freedom regarding their behaviour in time. 
Indeed if a scalar function $\alpha(x)$ satisfies 
\begin{equation}
\label{neu}
  \left\{
    \begin{aligned}
\Delta_x \,  \alpha & = 0  &&  \quad  \textrm{in } \Omega, 
\\  \partial_n   \,  \alpha  & = 0  &&  \quad \textrm{on }  \BorderNav ,
  \end{aligned}
    \right.
\end{equation}
then for any function $\eta(t)$, the vector field $\eta(t) \nabla_x  \alpha (x)$ satisfies the first three equations of \eqref{Eeq} for an appropriate pressure. 
In particular it is possible to choose a nonzero function  $\eta(t) $ satisfying $\eta(0) = \eta(T) = 0 $, so that this process leads to a field $u^0$ starting from zero initial data and  which vanishes at time $T$.   
Moreover the set of the scalar functions satisfying the underdetermined Neumann problem \eqref{neu} is rich enough to provide, by an appropriate gluing strategy, vector fields flushing the whole domain on the time interval $[0,T]$. 
\begin{lem} \label{lemma.euler}
	There exists a smooth solution~$(u^0, p^0)$ to 
\begin{gather}
  \label{Eeq0}
    \left\{
    \begin{aligned}
     \partial_t u^0 + \left( u^0 \cdot \nabla \right) u^0 + \nabla p^0 & =   0  &&\quad  \textrm{in } \Omega,
\\ \diverg u^0  & =  0 &&\quad  \textrm{in } \Omega,
\\   u^0 \cdot n  & =  0  &&\quad \textrm{on }  \BorderNav ,
\\   u^0 \rvert_{t = 0}  & = 0 &&\quad   \textrm{in}     \Omega ,
\\   u^0 \rvert_{t = T}  & = 0 &&\quad   \textrm{in} \Omega ,
\end{aligned}
    \right.
\end{gather}
such that the smooth solution $(u^1, p^1)$ to 
system \eqref{linu} satisfies  $u^1 \rvert_{t = T}  = 0$  in $ \Omega$. 
\end{lem}

This lemma is proved on the one hand  by the first author in  the 
papers~\cite{MR1233425}
and~\cite{MR1380673} respectively  for 2D simply connected domains and for general 2D domains 
{when} $\Sigma$ intersects all connected components of $\partial\Omega$,  and on the other hand by Glass in 
 \cite{MR1485616} and \cite{MR1745685} for the corresponding cases in 3D.

In the sequel, when we need it, we will implicitly extend the previous fields $u^0$ and $u^1$ by zero after $T$.

\section{Boundary layer}

The difficulty comes from the fact that the Euler equation, which models 
the behavior of a perfect fluid, not subject to friction, is only associated 
with the  boundary condition  $u\cdot n = 0$ for an impermeable wall
 and does not satisfy in general the Navier slip-with-friction boundary condition on $\BorderNav$.
An accurate description of  a solution of the Navier-Stokes equation near $\BorderNav$, even for a small viscosity, has therefore to use an expansion where a corrector is added to a solution to the Euler equation. 

\medskip
In the uncontrolled setting for which $\Sigma = \emptyset$ 
the description of the  behavior of the Navier-Stokes equation under the Navier
slip-with-friction condition in the  vanishing viscosity limit was performed by  Iftimie and the third author 
in  \cite{MR2754340} thanks to a multiscale  asymptotic expansion  involving a boundary layer term $v$ of 
amplitude~$\mathcal{O}(\sqrt{\varepsilon})$ and of thickness~$\mathcal{O}(\sqrt{\varepsilon})$
for a vanishing viscosity $\varepsilon$. 
Let us first briefly recall this result which will be extended in the sequel to the controlled case. 

Let us use here temporarily again the notation $u^0$ for a smooth solution to the Euler equations on the time interval $[0,T]$ with the impermeability condition 
 $u^0 \cdot n  =  0$ on the full boundary $\partial\Omega$.  
 The boundary layer corrector will involve an extra variable describing  the fast variations of the fluid velocity in the normal  direction near the boundary  and will 
 be given as a solution to an initial boundary value problem with a boundary condition with respect to this extra variable. 
We introduce a smooth function $\varphi : \R^d \rightarrow \R$ such 
that $\varphi = 0$ on $\partial \Omega$, $\varphi > 0$ {in} $\Omega$ and $\varphi 
< 0$ outside of~$\bar{\Omega}$. Moreover, we assume that ${|\varphi(x)|} = 
\mathrm{dist}(x, \Omega)$ in a small neighborhood of $\partial \Omega$. Hence, 
the normal~$n$ can be computed as $- \nabla \varphi$ close to the 
boundary and extended smoothly within the full domain~$\Omega$. 
   The notation $\tanpart{\cdot}$, introduced in \eqref{not-tan},  is extended  accordingly. 
We also introduce the following definitions:
\begin{equation*} 
 u^0_\flat(t,x)  {:=} - \frac{u^0(t,x) \cdot n(x)}{\varphi(x)}  \,  \text{ and  }  \,  g^0(t,x)  {:=} 2 \chi(x) N(u^0)(t,x) \quad    \textrm{ for } x  \in \Omega,
 \end{equation*} 
where   $\chi$ is a smooth {cut-off} function  satisfying $\chi = 1$ on $\partial \Omega$.
 Even though $\varphi$ vanishes on~${\partial \Omega}$, $u^0_\flat$ is 
    not    singular near the boundary because of the impermeability condition 
    $u^0\cdot n= 0$. Indeed since $u^0$ is smooth, a Taylor expansion proves
    that $u^0_\flat$ is smooth in~$\bar{\Omega}$.
The boundary layer corrector will be described by a smooth vector field $v$ 
 expressed in terms both 
of the slow space variable $x \in \Omega$ and a fast scalar variable 
$z = \varphi(x)/\sqrt{\varepsilon}$, where $v(t,x,z)$
satisfies the equation:
\begin{equation}
 \label{eq.v}
  \partial_t v + \tanpart{(u^0 \cdot \nabla) v + (v \cdot \nabla) u^0}  + u^0_\flat z \partial_z v - \partial_{zz} v  =  0,
\end{equation}
for $x$ in $\bar{\Omega}$ and $z$ in $ \R_+$, with the following boundary condition at $z=0$: 
\begin{equation}
 \label{eq.bv}
   \partial_z v(t, x, 0)  =  g^0(t,x) . 
   \end{equation}
We refer to  \cite[Section 2]{MR2754340} for a detailed heuristic  of the equations  \eqref{eq.v} and \eqref{eq.bv}.
Let us only mention here that these equations are obtained by plugging 
\begin{equation*} 
 u^0(t,x)  + \sqrt{\varepsilon} v\left(t,x, \frac{\varphi(x)}{\sqrt{\varepsilon}}\right)    \text{ instead of } \ue(t,x) 
 \end{equation*}
 into  the first and fourth equations of  \eqref{niou} and keeping the terms of higher order (taking into account that $u^0$ satisfies the Euler equations). 
Indeed the pressure  $ p^\varepsilon $ has to be expanded as well, into the sum of the Euler pressure and of a boundary layer term 
but the latter can be eliminated from the resulting equation
by distinguishing the normal and tangential parts. Thus this pressure  boundary layer  term acts as  a projection on the convective terms and this is  why  the second term in  \eqref{eq.v} is only tangential. 

The Cauchy problem associated with \eqref{eq.v} and \eqref{eq.bv}
   is well-posed in Sobolev spaces. Moreover for any $x\in \bar{\Omega}$, $z \geq 0$ and $t\geq 0$,
    we have 
    \begin{equation}
  \label{orto}
    v(t,x,z) \cdot n(x) = 0
\end{equation}
     It is easy to check that the solution inherits this condition 
    from the initial and boundary data. 
    This orthogonality property is the reason why equation~\eqref{eq.v} is linear.
    Indeed, the quadratic term $(v\cdot n) \partial_z v$ should have been taken
    into account if it did not vanish.
   Thanks to the cut-off function $\chi$, satisfying $\chi = 1$ on $\partial \Omega$, $v$ is compactly supported in $x$ near $\partial \Omega$, while ensuring
    that $v$ compensates the Navier slip-with-friction boundary trace of~$u^0$.

Then it is proved  in \cite{MR2754340} that the Leray solutions $ \ue$ to the  Navier-Stokes equation  can be described by  the following expansion in $L^\infty \big((0,T); L^2 (\Omega)\big)$:
\begin{equation*} 
 \ue(t,x) 
  = u^0(t,x) 
 + \sqrt{\varepsilon} v\left(t,x, \frac{\varphi(x)}{\sqrt{\varepsilon}}\right) +  O(\varepsilon) .
 \end{equation*}
Let us highlight that this expansion holds up to any time $T > 0$ for which 
$u^0$ is a smooth solution to the Euler equations on the time interval $[0,T]$.
On the other hand this analysis fails to describe the vanishing viscosity limit of the  Navier-Stokes equation 
for large times of order $O(\frac1{\varepsilon})$, even in the case where the Euler solution stays smooth for all times.

\medskip
Now going back to the controlled setting for which $\Sigma \neq \emptyset$ we expect to be able to describe 
 the  behavior of the Navier-Stokes equation near the uncontrolled part $\BorderNav$ of the boundary 
 in the  vanishing viscosity limit thanks to a similar expansion. 
 Indeed since we aim at finding a Navier-Stokes solution satisfying  \eqref{newgoal} we consider 
the following refined expansion:
\begin{align} 
 \label{eq.expansion}
 \ue(t,x) 
 & = u^0(t,x) 
 + \sqrt{\varepsilon} v\left(t,x, \frac{\varphi(x)}{\sqrt{\varepsilon}}\right)
 + \varepsilon u^1(t,x) 
 + \varepsilon \Rem(t,x), 
\end{align}
where  $u^0$ and $u^1$  are as  Lemma \ref{lemma.euler} and the vector field   $\Rem$ is wished to be $o(1)$ at time $\frac{T}{\varepsilon}$.
If so, and since the fields $ u^0$ and $ u^1$ are zero  after the time $T$, the leading part of the 
the expansion  \eqref{eq.expansion}   after $T$  is given by the second term in the right hand side and we therefore must understand
the large time behavior of this boundary layer.
For $t \geq T$, the equations \eqref{eq.v} and \eqref{eq.bv} reduce to 
\begin{equation} \label{eq.v.after}
  \left\{
  \begin{aligned}
    \partial_t v - \partial_{zz} v  & = 0,
    && \quad \textrm{ for } z \in \R_+ , \\
    \partial_z v(t, x, 0) & = 0 , 
  \end{aligned}
  \right.
\end{equation}
where the slow variables $x \in \Omext$ play the role 
of parameters through the  ``initial''  data $ \bar{v}(x,z) :=  v(T,x,z)$. 

This heat system  dissipates 
towards the null equilibrium state.  
Unfortunately the natural decay 
(that is without any 
assumption on $ \bar{v}$) 
at the final time $t = T/\varepsilon$  only yields 
\begin{equation} \label{eq.relax.v.naive}
 \left\| \sqrt{\varepsilon} 
 v\left(\frac{T}{\varepsilon},\cdot,\frac{\varphi(\cdot)}
 {\sqrt{\varepsilon}}\right) \right\|_{L^2(\Omega)}
 = \mathcal{O}\left(\varepsilon \right) ,
\end{equation}
which is  not sufficient in view of  the wished estimate   \eqref{newgoal}.
Physically, this is due to the fact that the average of $v$ is preserved under 
its evolution by equation~\eqref{eq.v.after} and to the fact that the energy contained by low 
frequency modes decays slowly.

\section{Well-prepared dissipation method}
\label{wpdm}

In order to overcome the previous difficulty we are going to  use the  \textit{well-prepared dissipation method}  first introduced 
in~\cite{MR3227326} by the second author in order to obtain a new  controllability result of the 1D  Burgers equation in the presence of a boundary layer. 
We will here adjust the method to the boundary layers associated with the Navier conditions in the vanishing viscosity limit of the Navier-Stokes equations.  
The idea is to design a  control strategy in order to enhance the natural dissipation of the boundary layer   after the time $T$.
Our strategy will be to guarantee that $\bar{v}$ satisfies a finite number of 
vanishing moment conditions for $k \in \N$ of the form:
\begin{equation} \label{eq.moments.k}
 \forall x \in \Omega, \quad \int_{\R_+} z^k \bar{v}(x,z) \dz = 0.
\end{equation}
This will allow to enhance the dissipation and to improve the estimate \eqref{eq.relax.v.naive} into 
\begin{equation} \label{eq.relax.v.na}
 \left\| \sqrt{\varepsilon} 
 v\left(\frac{T}{\varepsilon},\cdot,\frac{\varphi(\cdot)}
 {\sqrt{\varepsilon}}\right) \right\|_{L^2(\Omega)}
 = o\left(\varepsilon \right) .
\end{equation}
 Actually we aim at constructing the different fields mentioned so far by restriction to the physical domain $\Omega$ of solutions to analogous problems in an larger domain $\mathcal O$  extended across $\Sigma$ (and $\mathcal O$ can be chosen smooth, bounded and simply connected) with source terms compactly supported in the added portion of the domain. 
This means in particular that we intend to find a solution that we still denote  $(\ue,p^\varepsilon)$ of the following Navier-Stokes equations: 
$$ \partial_t \ue + \left( \ue \cdot \nabla \right) \ue - \varepsilon \Delta \ue + \nabla p^\varepsilon   =  \zeta^\varepsilon $$
for $x$ in $\mathcal O$ where the source term $ \zeta^\varepsilon (t,x)$ is a vector field supported for $x$ in  $\mathcal O \setminus \bar{\Omega}$, of the form 
\begin{align*} 
 \label{eq.expansion.zeta}
  \zeta^\varepsilon
 & =  \zeta^0(t,x) 
 + \sqrt{\varepsilon}  \zeta_v \left(t,x, \frac{\varphi(x)}{\sqrt{\varepsilon}}\right)
 + \varepsilon \zeta^1(t,x) , 
\end{align*}
where $ \zeta^0$ and $ \zeta^1$ are smooth vector fields used in order to insure Lemma \ref{lemma.euler} whereas the vector field 
$  \zeta_v (t,x,z)$ is devoted to the control of the moments of the boundary layer. 
Indeed we now aim at obtaining a profile $v$ solution to the following equation: 
\begin{equation}
 \label{eq.vE}
  \partial_t v + \tanpart{(u^0 \cdot \nabla) v + (v \cdot \nabla) u^0}  + u^0_\flat z \partial_z v - \partial_{zz} v  = \zeta_v   ,
\end{equation}
for $x$ in $\mathcal O$ and $z$ in $ \R_+$. 
Since the initial boundary value satisfied by $v$ is linear, its moments at time $T$ (see the left hand side of  \eqref{eq.moments.k}), can be decomposed as the sum of an addend due to the right hand side of  \eqref{eq.bv} and of an addend due to the outside control (see the right hand side of  \eqref{eq.vE}), which generates some moments outside, and are convected inside the domain by the field $u^0$, see the second term in  \eqref{eq.v}.
 Indeed, according to Duhamel's formula, the second addend is given by an integral over the time interval $[0,T]$, which allows to insure the  condition  \eqref{eq.moments.k} for all $x$ in $\Omega$.

 \begin{figure}[ht!]
     \begin{center}
         \includegraphics{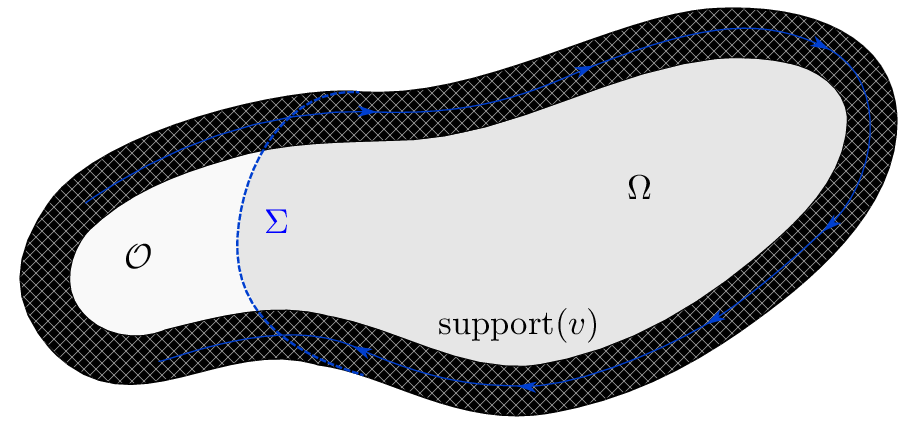}
        \end{center}
        \caption{\label{Figure:kaitenzushi}The kaitenzushi strategy}
    \end{figure}

Let us use here the following metaphor: see the extended domain as a conveyor-belt sushi restaurant, the added part of the extended domain as the kitchen and the moments as the plates (see Figure~\ref{Figure:kaitenzushi}). In order to send some plates from the kitchen, without sending the chef into the dining room, we use the transport by the field $u^0$ as a conveyor belt to serve the wished moments (compensating what comes from the uncontrolled part of the  boundary) all along the boundary before the end of the service, which corresponds here to the time  $T$. In this process it is crucial to maintain the orthogonality condition \eqref{orto} (otherwise the linearity of the equation would dramatically fall down because of the term $(v\cdot n) \partial_z v$ mentioned above). 
Let us observe that it seems impossible to control completely the boundary layer $v$ because the plates need some time to be conveyed from the kitchen and are therefore strongly regularized in $z$ (the equation \eqref{eq.bv} is parabolic in $z$) when they are supposed to be compensating what comes instantly from the nonhomogeneous data on the uncontrolled part of the boundary and is therefore far less regularized. Thus a compensation is only possible for a projection on a  functional space containing the two types of contributions  and we precisely make use of some finite dimensions projections by adjusting a finite number of moments.

\section{Estimates of the remainder}

Going back to the  velocity expansion  \eqref{eq.expansion} we are led now to the issue of estimating the large time behaviour of the remainder $ \Rem$.  
The field   $ \Rem$ can be naturally defined as the solution to a Navier-Stokes type equation of the form: 
\begin{equation} 
\label{eqR}  \partial_t \Rem 
   + \left( \ue \cdot \nabla \right) \Rem
   - \varepsilon \Delta \Rem
   + \nabla \pi^\varepsilon 
   + A^\varepsilon \Rem
    = {f^\varepsilon}  , 
  \end{equation}
  where    $\pi^\varepsilon$ denotes the pressure associated with  the vector field   $\Rem$, the notation $A^\varepsilon$ stands for an amplification operator and 
 $f^\varepsilon$ for a source term both due to the terms which were omitted in the equations of 
$u^0$, $u^1$ and $v$ for being of higher order in $\varepsilon$.  
 Since the field $ u^1$ bears the initial data $u_*$, this remainder starts with a zero initial data (taking the trace at the initial time of the equality \eqref{eq.expansion} and taking into account that $u^0$ and $v$ start with zero initial data) but is generated by the source term $f^\varepsilon$ and possibly amplified by mean of the term $A^\varepsilon \Rem$. 
 Of course the equation \eqref{eqR} is completed with the divergence free condition and  some initial and boundary conditions. 
 Here the key points in the large time estimate of $\Rem $ 
 are on the one hand that the quadratic nonlinearity in term of $\Rem$ which is hidden in the third term of  \eqref{eqR}  is tamed by a factor $\varepsilon$, see the   velocity expansion  \eqref{eq.expansion},  and on the other hand that the effects of both $A^\varepsilon$ and  $f^\varepsilon$ are tamed by the enhanced dissipation hinted in Section  \ref{wpdm}. 
Let us refer once more to  \cite{CMS} for more on the technicalities and only conclude here that the result of an energy estimate is that 
\begin{equation}
  \label{newgoalR}
 \left\| \Rem \left( \frac{T}{\varepsilon} ,\cdot\right)\right\|_{L^2(\Omega)}=o(1) .
\end{equation}
%

\section{Conclusion}

\begin{figure}[ht!]
    \begin{center}
        \includegraphics{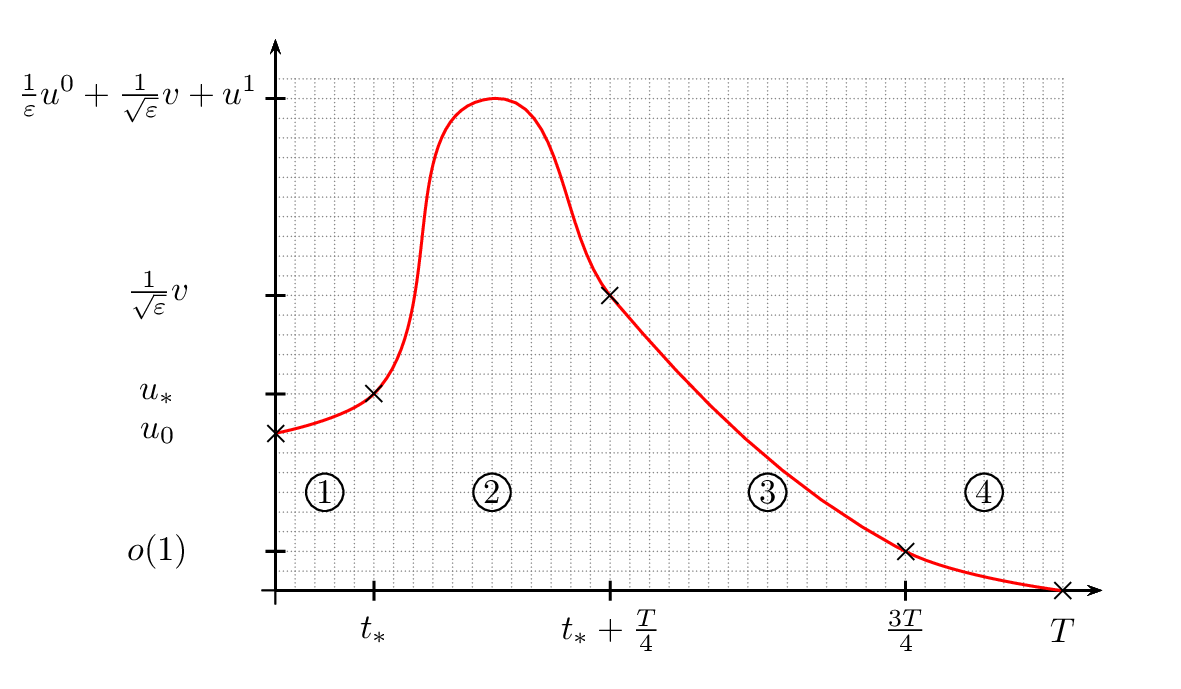}
    \end{center}
    \caption{\label{Figure:graph2} Four main steps of the evolution}
\end{figure}

Taking into account that the fields $u^0$ and $u^1$ vanish after $T$,  estimates \eqref{eq.relax.v.na} and~\eqref{newgoalR} plugged into expansion  \eqref{eq.expansion} yield  \eqref{newgoal} and therefore conclude the proof of Theorem  \ref{main} thanks to the preliminary reductions performed in Sections~\ref{sec-red} and~\ref{sec-ff}. 
The evolution of the state during the control strategy is 
pictured in Figure~\ref{Figure:graph2}.

\noindent The main steps of the proof can be summarized as in Table~\ref{Table}.

\begin{table}[ht]
\centering 
\begin{tabular}{c c c c c} 
\hline\hline 
Stage &  References  &  Active control & Linear behaviour &  Used effect\\ [0.5ex] 
\hline 
1 &  \cite{Leray}  & No &  & Dissipation \\ 
2 & \cite{MR1233425,MR1380673,CMS,MR1485616,MR1745685} & Yes &  $v$ & Convection \\
3 & \cite{CMS, MR3227326}   & No &  $v$,  $r^\varepsilon$  & Dissipation \\
4 & \cite{MR2224824, MR2268275} & Yes &  $u^\varepsilon$ & Dissipation \\ [1ex] 
\hline 
\end{tabular}
\caption{\label{Table}Main features of the control steps} 
\end{table}

\section{Perspectives}

Let us mention a few questions inspired by this work:
\begin{enumerate}
	\item Provided a smooth initial data, is there a strong solution to the 3D Navier-Stokes system reaching zero at time $T$?
	 The 2D case follows from Theorem~\ref{main}.
	\item Is it possible to deduce from the previous analysis some Lagrangian controllability results ? This would extend the results obtained  in~\cite{MR2579376},~\cite{MR3022084} 
	 for the incompressible Euler equations and in ~\cite{2016arXiv160203045G} for the stationary Stokes 
	equation.This issue is actually related to the previous one as  Lagrangian setting requires enough regularity for the flow to be  controlled. 
		\item Last but not least, is it possible to tackle \hyperref[lbOP]{(OP)} in the more difficult 
case of the no-slip  boundary condition, at least for some favorable geometric 
settings? This is a very challenging open problem because the no-slip  boundary 
condition gives rise to boundary layers that have a larger amplitude than 
Navier slip-with-friction boundary layers. 
We refer to the nice recent survey  \cite{2016arXiv161005372M}
by  {Maekawa} and {Mazzucato} for more on boundary layers in the no-slip case.
\end{enumerate}


\subsection*{Acknowledgements}

The two first authors 
 were partly supported by  ERC Advanced Grant 266907 (CPDENL) of  the 7th Research Framework Programme (FP7). 
The third author  thanks  the Agence Nationale de la Recherche, Project DYFICOLTI, grant ANR-13-BS01-0003-01, Project IFSMACS, grant ANR-15-CE40-0010  for their financial support and the third author thanks the hospitality of RIMS during the workshop on ``Mathematical Analysis of Viscous Incompressible Fluid".

\newpage

Jean-Michel \textsc{Coron},
Sorbonne Universit\'es, UPMC Univ Paris 
06, Lab. Jacques-Louis Lions, UMR CNRS 7598; 
\email{coron@ann.jussieu.fr}, 

\bigskip
Fr\'ed\'eric \textsc{Marbach},
Sorbonne Universit\'es, UPMC Univ Paris 
06, Lab. Jacques-Louis Lions, UMR CNRS 7598;  
\email{marbach@ann.jussieu.fr},

\bigskip
Franck \textsc{Sueur},
Institut de Math\'ematiques de Bordeaux, 
UMR CNRS 5251, Universit\'e de Bordeaux; 
\email{Franck.Sueur@math.u-bordeaux.fr}.

\end{document}